\documentclass[12pt]{amsart}
\usepackage{amssymb}
\usepackage{latexsym}
\def\R{\mathbb R}
\def\C{\mathbb C}

\def\N{\mathbb N}

\def\OO{{O}}

\def\re{\operatorname{Re}}

\def\arg{\operatorname{arg}}

\newtheorem{la}{Lemma}

\newtheorem*{thm}{Theorem}

\theoremstyle{remark}
\newtheorem*{ack}{Acknowledgments}
\newtheorem*{rem}{Remark}
\begin{document}
\title{Baker domains for Newton's method}
\author{Walter Bergweiler}
\thanks{The first author was supported by the Alexander von
Humboldt Foundation and by the G.I.F.,
the German--Israeli Foundation for Scientific Research and
Development, Grant G -809-234.6/2003.}
\address{Mathematisches Seminar,
Christian--Albrechts--Universit\"at zu Kiel,
Ludewig--Meyn--Str.\ 4,
D--24098 Kiel,
Germany}
\email{bergweiler@math.uni-kiel.de}
\author{D.\ Drasin}
\address{Department of Mathematics,
Purdue University, West Lafayette, IN 47907, USA}
\email{drasin@math.purdue.edu}
\author{J.~K.~Langley}
\address{School of Mathematical Sciences,
University of Nottingham, NG7 2RD, UK}
\email{jkl@maths.nott.ac.uk}
\begin{abstract}
We show that there exists an entire function without 
finite asymptotic
values for which the associated  Newton function tends to infinity
in some invariant domain.  The question whether such a function
exists had been raised by Douady.
\end{abstract}
\maketitle

\section{Introduction and result}
Let $f$ be an entire function. Newton's method
for finding the zeros of $f$ consists of
iterating the function
$$N(z):=z-\frac{f(z)}{f'(z)}.$$
If $\xi$ is a zero of $f$, then $N(\xi)=\xi$ and $|N'(\xi)|<1$,
so there is an
$N$-invariant domain $U$ containing $\xi$
in which the iterates $N^k$ of $N$ converge to $\xi$
as $k\to\infty$.
(Here $N$-invariance of $U$ means that $N(U) \subset U$.)

There may also be $N$-invariant domains
in which the iterates of $N$ tend to $\infty$. A simple example
is given by $f(z)=P(z) \exp Q(z)$ where $P$ and $Q$ are polynomials,
with $Q$ nonconstant. Then $N$ is rational. Moreover, in the
terminology of complex dynamics, the point at $\infty$ is 
a fixed point of multiplier $1$ of $N$, and the iterates of $N$ 
tend to $\infty$ in the Leau petals associated to this fixed point.

If $f$ does not have the above form, then $N$ is transcendental;
see~\cite{Ber93} for an introduction to the iteration theory
of transcendental meromorphic functions.
A maximal $N$-invariant domain where the iterates of $N$ tend
to $\infty$ is called an {\em invariant Baker domain}.

A simple example (cf.~\cite{BHKMT}) is given by functions
$f$ for which $f(z)\sim\exp(-z^n)$ as $z\to\infty$ in some 
sector $|\arg z|<\varepsilon$. Then 
$N(z)=z+(1/n+o(1))z^{1-n}$
and this implies that $N^k|_U\to\infty$ 
as $k\to\infty$
for some $N$-invariant domain $U$ containing all sufficiently large real
numbers. Note that $f(x)\to 0$ as $x\to\infty$, $x\in\R$.
Thus $0$ is an asymptotic value of $f$, the positive real axis being
an asymptotic path. 
Figuratively speaking one might say that 
Newton's method believes that there is a zero of $f$ at $+\infty$,
and thus it yields a domain $U$ containing all sufficiently large 
reals such that $N^k(z)\to+\infty$ for $z\in U$ as $k\to\infty$. 

The question arises 
whether an entire function $f$ must always have $0$ as an asymptotic
value if $N$ has an invariant Baker domain. 
This question was raised by A.\ Douady and has been brought to our
attention by J.\ R\"uckert.
It has been shown by
X.\ Buff and J.\ R\"uckert~\cite{BR}
that the answer to this question is positive
in situations much more
general than those given above. However, we shall show
that this is not always the case.

\begin{thm} 
There exists an entire function $f$ without finite 
asymptotic values such that $N(z)=z-f(z)/f'(z)$ has
an invariant Baker domain.

Moreover, $f$ can be chosen to be of any preassigned
order strictly between $\frac12$ and~$1$.
\end{thm}

We explain the basic idea of the construction.
Using functions of the type introduced by S.\ K.\ Bala\v{s}ov~\cite{Bal},
in \S 2 we construct an 
entire function $f$ 
of order less than $1$ (and in fact of any preassigned
order strictly between $\frac12$ and~$1$)
which satisfies 
$$f(z)\sim \sqrt[q]{z}$$
for some integer $q$ and some branch of the $q$-th root 
as $z\to\infty$ in
the spiraling region
$$S:=\left\{ re^{ic\log r +i\theta}:
r>1,
|\theta|<\theta_0\right\},$$
where $c:=\pi/\log (q-1)$ and $0<\theta_0<\pi$. 
Here the relation between $c$ and $q$ is such that
$S$ is invariant under $z\mapsto -pz$
where $p:=q-1$.
We shall see in~\S 3 
that
$$\frac{f'(z)}{f(z)}\sim 
\frac{1}{qz}$$
so that
$$N(z)=z-\frac{f(z)}{f'(z)}\sim
-pz.$$
This will show that $S$ contains an $N$-invariant domain in which the
iterates of $N$ tend to $\infty$.
Hence $N$ has an invariant Baker domain. Finally we shall show
in~\S 4, using the Denjoy-Carleman-Ahlfors Theorem, that $f$ has
no finite asymptotic values.

\section{The construction of $f$}\label{construction}
Let $(a_k)$ be a sequence of complex numbers tending 
to infinity. For $r>0$ let $n(r)$ be the 
number of $a_k$, taking account 
of repetition, in $|z| \leq r$.
Let
$$\rho:=\limsup_{r\to\infty}\frac{\log n(r)}{\log r}.$$
Equivalently, $\rho$ is the exponent of convergence of the
sequence $(a_k)$. It is well known that the canonical product
$\Pi$ whose zeros are the $a_k$ has order $\rho$; that is,
$$\rho=\limsup_{r\to\infty}\frac{\log \log M(r,\Pi)}{\log r}$$
where $M(r,\Pi):=\max_{|z|=r}|\Pi(z)|$ is the maximum modulus.
There are standard results concerning the asymptotic behavior
of $\Pi$ if all $a_k$ lie on one ray and if
\begin{equation} \label{2a}
n(r)\sim \Delta r^\rho
\end{equation} 
for some $\Delta>0$ as $r\to\infty$.
These results have been extended by S.\ K.\ Bala\v{s}ov~\cite{Bal}
to the case where the $a_k$ lie on a logarithmic spiral,
say
\begin{equation} \label{2b}
a_k\in\left\{re^{ic\log r}: r\geq 1 \right\},
\end{equation} 
where $c>0$. 
We quote only a simplified version of
Bala\v{s}ov's result~\cite[Theorem~1]{Bal},
as this suffices for our purposes.
\begin{la} \label{balasov}
Let $(a_k)$ be a sequence satisfying~$(\ref{2a})$ and~$(\ref{2b})$.
Suppose that $\rho$ is not an integer.
Let $\Pi$ be the canonical product formed  with the $a_k$. Then
$$\lim_{r\to\infty}\frac{\log \Pi(re^{ic\log r+i\theta})}{r^\rho}
=
-\frac{2\pi i \Delta \exp\left(i\rho\theta/(1+ic)\right)}
{1-\exp\left(i2\pi\rho/(1+ic)\right)}
$$
for $0<\theta<2\pi$ and
a suitable branch of the logarithm, the convergence being uniform
for $\varepsilon\leq \theta\leq 2\pi-\varepsilon$ if $\varepsilon>0$.
In particular, 
\begin{equation} \label{2c}
\lim_{r\to\infty}\frac{\log |\Pi(re^{ic\log r+i\theta})|}{r^\rho}
=h(\theta):=
-2\pi  \Delta\re\left(
\frac{i \exp\left(i\rho\theta/(1+ic)\right)}
{1-\exp\left(i2\pi\rho/(1+ic)\right)}
\right).
\end{equation} 
\end{la}
Now let $\frac12 < \rho<1$ and $\Delta>0$. Choose $p\in \N$ such that 
\begin{equation} \label{defmu}
\mu:=\frac{\rho}{1+c^2}:=\frac{\rho}{1+(\pi/\log p)^2}>\frac12,
\end{equation}
thus defining $c:=\pi/\log p$.
Note that since $\frac12 <\mu< \rho<1$ we have $c<1$ and hence
$p>\exp(\pi)>23$.
Let $(a_k)$ be a sequence satisfying~(\ref{2a}) and~(\ref{2b})
and let $\Pi$ be the canonical product formed with
the $a_k$ so that~(\ref{2c})
holds.

A series of elementary modifications
of $\Pi$ will produce the function $f$ of our theorem.

A computation shows that 
$$h(0)=-2\pi\Delta\re\left(
\frac{ i}{1- \exp\left(i2\pi\rho/(1+ic)\right)}
\right)
=\frac{2\pi\Delta\exp(2\pi\mu c) }{\left|1-
 \exp\left(i2\pi\rho/(1+ic)\right)\right|^2}
\sin(2\pi \mu).$$
Since $\frac12 < \mu<1$ we thus have $h(0)<0$.
Hence there exists $\theta_0>0$ such that $h(\theta)<0$ for 
$|\theta|<\theta_0$.
For 
$0<\varepsilon<\theta_1<\theta_0$ we thus deduce from~(\ref{2c})
that there exists $\eta_0>0$ such that
\begin{equation} \label{2d}
\log |\Pi(re^{ic\log r+i\theta})|\leq -\eta_0r^\rho \quad
\text{for }\  \varepsilon\leq |\theta|\leq\theta_1,
\end{equation} 
provided $r$ is sufficiently large.

We show that an estimate of this type also
holds for $|\theta|<\varepsilon$.
In order to do so,
we use a standard estimate which in slightly different form
can be found
in~\cite[p.~548]{Hay89} or~\cite[p. 117]{Tsu}.
\begin{la} \label{latsuji}
Let $D\subset \C$ be an unbounded domain. For $r>0$ such
that the circle $C_r:=\{z\in\C:|z|=r\}$ intersects $D$, let
$r\theta(r)$ be the linear measure of that intersection.
Let $\theta^*(r):=\theta(r)$ if $C_r\not \subset D$ and
let $\theta^*(r):=\infty$ and thus $1/\theta^*(r):=0$ if 
$C_r\subset D$. 

Suppose that $u:\overline{D}\to [-\infty,\infty)$
is continuous in $\overline{D}$ and subharmonic in $D$.
Suppose also that $u$ is bounded above
on $\partial D$, but not bounded above in $D$.
Let $0<\kappa<1$
and let $R>0$ be such that $C_R$ intersects $D$.
Then 
$B(r,u):=\max_{|z|=r} u(z)$ satisfies
$$\log
B(r,u)
\geq\pi\int_R^{\kappa r}\frac{dt}{t\theta^*(t)}-\OO(1)$$
as $r\to\infty$.
\end{la} 
We may assume that 
$\varepsilon$ in~(\ref{2d})
is chosen such that $0<\varepsilon<\pi/2$. We consider the 
spiralling domain
$$D:=\left\{ re^{ic\log r +i\theta}:
r>1,
|\theta|<\varepsilon\right\}$$
and the function 
$$
u(z):=\log |\Pi(z)|+\eta_0|z|^\rho.$$
Then $u$ is continuous in $\overline{D}$, subharmonic in $D$
and bounded above on $\partial D$. 
We claim that $u$ is also bounded above in $D$.  Otherwise, on applying
Lemma~\ref{latsuji}
and noting that $\theta^*(r)=2\varepsilon$
we find that
$$\log
B(r,u)
\geq\pi\int_R^{\kappa r}\frac{dt}{2\varepsilon t}-\OO(1)
=\frac{\pi}{2\varepsilon}
\log r
-\OO(1)
>\log r$$
and thus 
$$\log M(r,\Pi)
=B(r,u)-\eta_0r^\rho>
r-\eta_0r^\rho>\frac{r}{2}$$
for large $r$. This implies that the order of $\Pi$ is at least~$1$,
a contradiction.

Thus $u$ must be bounded above in $D$, and 
this, together with~(\ref{2d}),  implies that 
if $0<\eta_1<\eta_0$, then
\begin{equation} \label{2e}
\log |\Pi(re^{ic\log r+i\theta})|\leq -\eta_1r^\rho \quad
\text{for }\   |\theta|\leq\theta_1
\end{equation} 
and sufficiently large~$r$.

Let $L$ be the natural parametrization of the logarithmic spiral
on which the $a_k$ lie; that is,
$L:[1,\infty)\to \C$, $L(t)=te^{ic\log t}$. Then
$\Pi(L(t))\to 0$ as $t\to\infty$ by~(\ref{2e}).
Thus there exists $t_0>1$
such that $|\Pi(L(t))|<|\Pi(L(t_0))|$ for $t>t_0$.

The function $f$ of our theorem will now be defined as follows.
We put $z_0:=L(t_0)$ and define $g_0(z):=\Pi(z+z_0)$.
Next we put $q:=p+1$ and $g_1(z):=g_0(z^q)$, and define
$\sigma:[0,\infty)\to \C$ by $\sigma(t)=\sqrt[q]{L(t_0+t)-z_0}$,
where the branch is chosen such that
$\arg \sigma(t)=q^{-1} \arg(L(t_0+t)-z_0)$ for $t>0$.
We then define 
\begin{equation} \label{defg2}
g_2(z):=\int_0^z g_1(\zeta)^n d\zeta
=\int_0^z \Pi(\zeta^q+z_0)^n d\zeta,
\end{equation}
where $n\in\N$. It will follow easily that
$$a:=\int_\sigma g_1(z)^n dz$$
is finite for all $n\in\N$, and using a
result of W.\ K.\ Hayman~\cite[Lemma~1]{Hay69}
we will see that $a\neq 0$ if $n$ is sufficiently large.
For such a value of $n$ we 
then define $g_3(z):=g_2(z)/az$ and note that $g_3$ is of the
form $g_3(z)=g_4(z^q)$ for some entire function $g_4$.
The function claimed in the statement of
the theorem is
$$f(z):=z g_4(z)^{q-1}.$$

To prove that $f$ 
has the desired properties,
we will determine the asymptotic behavior
of the $g_j$ and $f$ in spiralling regions similar to $D$.
We first note that 
it follows from~(\ref{2e})
that if $0<\eta_2<\eta_1$ and if 
$0<\theta_2<\theta_1$, then
$$
\log |g_0(re^{ic\log r+i\theta})|\leq -\eta_2r^\rho \quad
\text{for }\   |\theta|\leq\theta_2
$$
and sufficiently large~$r$.
This implies that if 
$|\theta|\leq\theta_2/q$ and if $r$ is sufficiently large, then 
$$\log |g_1(re^{ic\log r+i\theta})|
=
\log |g_0(r^qe^{ic\log (r^q)+iq\theta})|
\leq -\eta_2r^{q\rho}.$$
With 
$$S_1:=\left\{ re^{ic\log r +i\theta}:
r>1,
|\theta|<\frac{\theta_2}{q}\right\}$$
we thus find that 
$$|g_1(z)|\leq \exp\left(-\eta_2|z|^{q\rho}\right)$$
if $z\in S_1$ is sufficiently large.
Moreover, we have 
$\sigma(t)\in S_1$ for large $t$, since
\begin{eqnarray*}
\arg \sigma(t)
&=&
\frac{1}{q}\arg(L(t_0+t)-z_0)\\
&=&
\frac{1}{q}\arg L(t_0+t)+o(1)\\
&=&
\frac{c}{q}\log  |L(t_0+t)|+o(1)\\
&=&
c \log|\sigma(t)|+o(1).
\end{eqnarray*}
We deduce that the integral
defining $a$ converges for all $n\in\N$.
In order to show that $a\neq 0$ for large $n$, we shall use the following
result of W.\ K.\ Hayman~\cite[Lemma~1]{Hay69}.
\begin{la} \label{lahayman}
Let $\gamma$ be a Jordan arc in 
$\C$ which tends to $\infty$ in both directions and
let $g$ be holomorphic in a domain containing $\gamma$.
Suppose that
$\int_\gamma |g(z)||dz|<\infty$ and that
$|g(z)|\to 0$ as $z\to\infty$ on $\gamma$.
Suppose also that $|g(z)|\leq M$ for $z$ on  $\gamma$, 
with equality for a single point $z_1$ on $\gamma$ which satisfies
$g'(z_1)=0$. Suppose finally that $\gamma$ cannot be deformed 
in a neighborhood of $z_1$ into a curve on which $|g(z)|<M$.
Then 
$$\int_\gamma g(z)^n dz\neq 0 $$
for all sufficiently large integers~$n$.
\end{la}
We apply this lemma with $g:=g_1$, 
the curve
$\gamma$ parametrized as 
$\gamma:\R\to\C$, 
$$\gamma(t):=
\left\{
\begin{array}{ll}
\sigma(-t) & \mbox{if $t\leq 0$,}\\[1mm]
\sigma^*(t):= e^{2\pi i /q}\sigma(t) & \mbox{if $t>0$,}
\\
\end{array}
\right.
$$
and $z_1:=\gamma(0)=0$.
Since by the choice of $z_0$ we have
$$|g_1( \sigma^*(t) )|=|g_1( \sigma(t) )|= | \Pi (L(t_0+t)) | 
<| \Pi (L(t_0)) | = |g_1( \sigma(0) )|$$
for $t>0$, it follows that $|g_1(\gamma(t))|<|g_1(z_1)|$ for $t\neq 0$.
Moreover, $g_1'(z_1)=g_1'(0)=0$, and thus
the hypotheses of Lemma~\ref{lahayman} are satisfied.
Since
$$\int_\gamma g_1(z)^ndz
=
-\int_\sigma g_1(z)^ndz+
\int_{\sigma^*} g_1(z)^ndz
=\left(-1+e^{2\pi i /q}\right)
\int_\sigma g_1(z)^ndz$$
we conclude from Lemma~\ref{lahayman} that 
$a=\int_\sigma g_1(z)^ndz\neq 0$ for sufficiently large 
values of~$n$.

We deduce that
$g_2(\sigma(t))\to a$
as $t\to\infty$.
More generally, $g_2(z)\to a$ as $z\to\infty$ in $S_1$. In fact,
if $z\in S_1$ then
$$g_2(z)-a=\int_{\tau_z} g_1(\zeta)^n d\zeta$$
for any path $\tau_z$ joining $z$ to $\infty$ in $S_1$.
For large $z\in S_1$ and suitable $\tau_z$ we find that
$$|g_2(z)-a| \leq 
\int_{\tau_z} |g_1(\zeta)|^n |d\zeta|\leq
\int_{\tau_z} 
\exp\left(-n\eta_2|\zeta|^{q\rho}\right)
|d\zeta|\leq
\exp\left(-\eta_3|z|^{q\rho}\right)
$$
for some $\eta_3>0$. 
It follows 
that if $z\in S_1$ is sufficiently large, then 
$$\left|g_3(z)-\frac1z\right|= \frac{|g_2(z)-a|}{|az|}
\leq\frac{\exp\left(-\eta_3|z|^{q\rho}\right)}{|az|} 
\leq\exp\left(-\eta_3|z|^{q\rho}\right).$$

Now let
$$S_2:=\left\{ re^{ic\log r +i\theta}:
r>1,
|\theta|<\theta_2\right\}.$$
For $z\in S_2$ we have $\sqrt[q]{z}\in S_1$ for a suitable 
branch. For large $z\in S_2$ we thus find that
$$\left|g_4(z)-\frac{1}{\sqrt[q]{z}}\right|= 
\left|g_3(\sqrt[q]{z})-\frac{1}{\sqrt[q]{z}}\right|
\leq \exp\left(-\eta_3|z|^{\rho}\right);$$
i.~e.~if $z\in S_2$ is sufficiently large, then 
\begin{equation} \label{2g}
\left|f(z)-\sqrt[q]{z}\right|
\leq \exp\left(-\eta_4|z|^{\rho}\right)
\end{equation}
for some $\eta_4>0$ and
a suitable branch.
\begin{rem}
We have chosen $z_0$ and $n$ in the way described only
to ensure that $a\neq 0$. In a generic situation we could probably
define $g_2$ via~(\ref{defg2}) with
$z_0=0$ and $n=1$.
\end{rem}
\section{Newton's method for $f$}
We choose $\theta_3$ with $0<\theta_3<\theta_2$ and define
$$S_3:=\left\{ re^{ic\log r +i\theta}:
r>1,
|\theta|<\theta_3\right\}.$$
Then there exists $\delta>0$ such that if $z\in S_3$ is sufficiently large,
then the closed disk of radius $\delta|z|$ around $z$ is contained in $S_2$.
With $d(z):=f(z)-\sqrt[q]{z}$ we deduce from~(\ref{2g}) that if $z\in S_3$ 
is sufficiently large, then 
\begin{eqnarray*}
\left|f'(z)-\frac{\sqrt[q]{z}}{qz}\right|
&=&\left|d'(z)\right|\\
&=&
 \frac{1}{2\pi}\left|\int_{|\zeta-z|=\delta|z|} \frac{d(\zeta)}{(\zeta-z)^2}
d\zeta \right|\\
&\leq&
 \frac{1}{\delta|z|}\max_{|\zeta-z|=\delta|z|} \left|d(\zeta)\right|\\
&\leq& 
\frac{1}{\delta|z|}\exp\left(-\eta_4(1-\delta)^\rho |z|^{\rho}\right)\\
&\leq& \exp\left(-\eta_5|z|^{\rho}\right)
\end{eqnarray*}
for some $\eta_5>0$.
Combining this with~(\ref{2g}) we find that if $z\in S_3$ is sufficiently large,
then
$$\left|\frac{f(z)}{f'(z)}-qz\right|
\leq \exp\left(-\eta_6|z|^{\rho}\right)$$
where $\eta_6>0$.
Recalling that $q=p+1$ we deduce that
\begin{equation} \label{npz}
|N(z)+pz| =
\left|z-\frac{f(z)}{f'(z)}+pz\right|
=
\left|\frac{f(z)}{f'(z)}-qz\right|
\leq \exp\left(-\eta_6|z|^{\rho}\right)
\end{equation}
for large $z\in S_3$.
In particular, 
$$\big| |N(z)|-p|z|\big| \leq \exp\left(-\eta_6|z|^{\rho}\right)$$
which implies that 
\begin{equation} \label{3a}
\big| \log |N(z)|-\log (p|z|)\big| \leq \exp\left(-\eta_6|z|^{\rho}\right)
\end{equation}
for large $z\in S_3$.
Moreover,~(\ref{npz}) yields
\begin{equation} \label{3b}
\left| \arg N(z)-\arg (-pz)\right| \leq \exp\left(-\eta_6|z|^{\rho}\right)
\end{equation}
for large $z\in S_3$.
Recalling that $c$ was chosen such that $c\log p=\pi$ we deduce 
from~(\ref{3a}) and~(\ref{3b}) that
\begin{eqnarray*}
&     & \big| \arg N(z)-c\log |N(z)|\big|
\\
&\leq & \big| \arg N(z)-\arg (-pz)\big| +
\big|\arg (-pz)- c\log  (p|z|)\big| \\
& & +
\big|c \log (p|z|)-c \log |N(z)|\big|
\\
&\leq & \big|\arg (-pz)- c\log  (p|z|)\big| + 
(1+c)\exp\left(-\eta_6|z|^{\rho}\right)
\\
&= & \big|\arg z+\pi - c\log  p- c\log |z|\big| + 
(1+c)\exp\left(-\eta_6|z|^{\rho}\right)
\\
&= & 
\big|\arg z- c\log |z|\big| +
(1+c)\exp\left(-\eta_6|z|^{\rho}\right)\\
&\leq & 
\big|\arg z- c\log |z|\big| +
\frac{1}{2|z|}
\end{eqnarray*}
for large $z\in S_3$.
Since $p>23$
we deduce from~(\ref{npz}) that
$|N(z)|> 2|z|$ if $z\in S_3$
and if $|z|$ is sufficiently large, say $|z|>r_0>1$.
Combining this with the previous estimate we conclude that if
$$\big|\arg z- c\log |z|\big|
< \theta_3-\frac{1}{|z|},$$
then
$$\big| \arg N(z)-c\log |N(z)|\big|
< \theta_3-\frac{1}{|z|}+\frac{1}{2|z|}
= \theta_3-\frac{1}{2|z|}
< \theta_3-\frac{1}{|N(z)|}$$
if  $z\in S_3$
and if $|z|$ is large enough, say
$|z|> r_1>r_0$. 
This implies that
$$U:=\left\{ re^{ic\log r +i\theta}:
r>r_1,
|\theta|<\theta_3-\frac{1}{r}\right\}$$
is $N$-invariant.
Since $|N(z)|\geq 2|z|$ for $z\in U$, we have
$|N^k(z)|\geq 2^k|z|$ for $z\in U$ and $k\in\N$. Thus
$N^k|_U\to\infty$ as $k\to\infty$.
Hence $U$ is contained in an invariant Baker domain of $N$.

\section{Asymptotic values of $f$}
Suppose that $f$ has a finite asymptotic value,
say
$f(z)\to b\in \C$ as $z\to \infty$ along a curve $\Gamma$.
The function 
$$F(z):=\frac{f(z)^q}{z}$$
is entire since $f(0)=0$.
By (\ref{2g})
we have $F(z)\to 1$ as $z\to \infty$ along the logarithmic spiral $L$ while 
$F(z)\to 0$ as $z\to \infty$ along $\Gamma$.
Thus $F$ has two finite
asymptotic values. By the Denjoy-Carleman-Ahlfors Theorem 
(see~\cite[\S XI.4.5]{Nev53}),
$F$~has order at least $1$. On the other hand, $F$ has the same order as $f$,
which has been taken less than $1$. This is a contradiction.
\begin{rem}
Our method will produce examples $f$ of any preassigned 
non-integer order $\rho>1$, as well as examples with more
than one invariant Baker domain.
We only sketch the modifications that have to be made.

We again choose $\rho$ and $p$ such that~(\ref{defmu}) is satisfied.
The condition $\mu<1$ need not be satisfied,
and there may be several, say $\ell$, intervals where $h(\theta)$ 
is negative and corresponding spiraling regions $S_1,\dots,S_\ell$
where $\Pi(z)\to 0$ as $z\to\infty$.
It is not difficult to see that
$\ell$ can be any given positive number.
Let $L_j$ be a curve starting at $0$ which outside the unit circle
is a logarithmic spiral in $S_j$ and which inside the unit circle is
a straight line from $0$ to the corresponding point of the unit circle.
Deforming one of the curves $L_j$ if necessary we may assume that 
there exists $z_0\in \bigcup_{j=1}^\ell L_j$ such that
$|\Pi(z_0)|>|\Pi(z)|$ for all $z\in \bigcup_{j=1}^\ell L_j$.
Defining $g_2$ by~(\ref{defg2}) for some large $n$ 
and defining $f$ then as in~\S\ref{construction}
we find that $f(z)\sim c_j\sqrt[q]{z}$ for some $c_j\neq 0$ 
as $z\to\infty$ in $S_j$.
As before we deduce that
$N(z)\sim -pz$ as $z\to\infty$ in $S_j$, $j=1,\dots,\ell$.
We thus obtain an entire function $f$ for which $N$ has $\ell$
invariant Baker domains.
A difference occurs in the proof that $f$ does not have finite
asymptotic 
values. Here we cannot simply appeal to the classical 
Denjoy-Carleman-Ahlfors Theorem, but 
instead use that the function
$f$ constructed has only $\ell$
``tracts''; see~\cite{Hay89}.

Bala\v{s}ov's result takes a different form if $\rho$ is an integer, 
but it seems 
possible to treat this case along the same lines.
\end{rem}

\begin{ack}
We thank Johannes R\"uckert for bringing Douady's question to our
attention, and for interesting discussions about it. We also 
acknowledge helpful conversations with Alexandre Eremenko.
\end{ack}

\end{document}